\documentclass[12pt]{amsart}
\usepackage{amssymb,latexsym,amsmath,amsthm,enumitem,mathrsfs,geometry}
\usepackage{fullpage}
\usepackage{fancyhdr}
\usepackage{xcolor}
\usepackage{bbm}
\usepackage{stmaryrd}
\usepackage{mathrsfs}
\usepackage{hyperref}
\usepackage{lineno}
\usepackage{diagbox}
\usepackage{array}
\usepackage{tikz}
\usetikzlibrary{arrows}
\usetikzlibrary{positioning, arrows.meta, shapes.geometric, calc}
\usepackage{tikz-cd}
\usepackage{graphicx}
\usepackage{float}
\usepackage{comment}
\usepackage{mathtools}
\usepackage{cleveref}

\newtheorem*{theorem*}{Theorem}

\newtheorem{conjecture}{Conjecture}[section]

\newtheorem*{remark*}{Remark}

\theoremstyle{definition}
\newtheorem{question}{Question}[section]
\newtheorem{hypothesis}{Hypothesis}[section]
\newtheorem{goal}{Goal}[section]

\begin{document}

\numberwithin{equation}{section}

\title{Mathematics in the age of AI}

\author{Terence Tao}
\address{UCLA Department of Mathematics, Los Angeles, CA 90095-1555.}
\email{tao@math.ucla.edu}


\date{\today}

\begin{abstract} An essay, based on a public lecture delivered at the 2026 International Congress of Mathematicians, on how the mathematical community might respond to the arrival of artificial intelligence tools that are capable of performing research-level mathematical tasks.  Rather than debating the capabilities of such tools, we condition on the hypothesis that these capabilities will arrive, and examine instead a question that is orthogonal to it: what the goals and values of mathematical research actually are.  The problem-solving component of mathematics is used as a case study.
\end{abstract}

\maketitle

\section{A historical prologue}

For centuries, mathematics operated successfully on ``naive'' foundations.  Practicing mathematicians proved theorems about sets, numbers, and infinities without feeling any particular need to say precisely what these quantites actually were (or what a proof itself was, for that matter); such questions were largely delegated to philosophers, and the working mathematician was free to get on with the mathematics.  

But in the early twentieth century, discoveries such as Russell's paradox in 1901 \cite{russell} and the G\"odel incompleteness theorems in 1931 \cite{godel} forced mathematicians to critically re-examine assumptions about their subject that had previously been left implicit.  The axioms of naive set theory contradicted each other; and a formal system could not simultaneously be consistent, sufficiently expressive, and capable of proving its own consistency.

The resulting \emph{crisis in foundations}, lasting roughly from 1900 to 1930, was a genuinely turbulent period for the subject.  But the end product of that turbulence was extremely valuable: an explicit, rigorous, and standardized foundational framework, in which the objects of mathematics and the rules for reasoning about them are laid out in a form that can be inspected, taught, and --- as it turns out --- mechanized.  There is certainly scope for further improvement in this framework; and foundational research continues to this day.  But our current foundations have survived a century of strenuous testing, and they now provide a trusted environment in which mathematics can be conducted with a very high degree of confidence.

I believe that we are now entering a era of comparable turbulence in mathematics.  This time, though, what is being stress-tested is not our foundational framework for mathematical \emph{truth}, but rather the largely implicit framework of mathematical \emph{values} and \emph{practices}: what we consider a contribution to be, what we reward, what we regard as understood, and who --- or what --- we regard as having done the work.  I argue that it will become necessary to make these unwritten goals of mathematics much more explicit; but once we have thoroughly examined and codified them, our community will emerge stronger and more resilient than before.

\section{The motivating question}

This article is organized around a single question.

\begin{question}[Community Response Question]\label{crq} How should the mathematical\footnote{The impacts of AI of course extend far beyond mathematics, but that is far too vast a topic to address here.} community respond to the advent of modern AI technologies, and their real and/or claimed capabilities to perform mathematical tasks?
\end{question}

This is a question for the entire community, and I do not presume to have all of the answers to it; nobody does.  Nevertheless, I have some things to say about how one might go about answering it.

The first thing to say is that Question \ref{crq} is \textbf{not} a mathematical question.  It is a metamathematical one, and also a political, ethical, sociological, and cultural one; it cannot be settled by pure logical argument.  However, in what follows I will deliberately \emph{borrow} the precise and familiar language of mathematics --- conjectures, hypotheses, and the like --- in order to clarify the structure of the question.

\section{The first subquestion: AI capability}

The answer to Question \ref{crq} depends crucially on a subquestion about what AI tools will actually be able to do.  It is convenient to formulate this pseudomathematically, not as a single conjecture, but as a \emph{family} of conjectures indexed by a large number of free parameters.

\begin{conjecture}[AI Capability Conjecture, template form]\label{aicc} At \emph{some} point in the near future, \emph{some} AI tools will, at \emph{some} expense, and with \emph{some} level of human supervision, be able to accomplish \emph{some} research-level mathematical tasks in \emph{some} fields of mathematics, with \emph{some} non-trivial success rate, and at \emph{some} level of correctness and quality.
\end{conjecture}

Each occurrence of the word ``some'' above should be read as a placeholder (or a free parameter).  One can obtain a great many distinct conjectures depending on how one fills these placeholders in.  The finer distinctions between these formulations are important, but they are not the point of this article; I will make only the coarse distinction between ``weak'' and ``strong'' forms of Conjecture \ref{aicc}.

If even weak forms of the AI Capability Conjecture turn out to be false, then we could safely dismiss the current generation of AI tools as being of no long-term significance to mathematical research, and largely continue with business as usual.  If, on the other hand, the strongest forms of the conjecture are true, then it becomes very challenging to maintain our current culture and practices unchanged --- particularly if we continue to prioritize such goals as obtaining as many solutions to unsolved problems as possible.

It is therefore difficult to have a constructive discussion on Question \ref{crq} while the status of Conjecture \ref{aicc} remains under dispute.  Unsurprisingly, then, most of the public debate about AI and mathematics has concerned which versions of the AI Capability Conjecture are true; I have myself devoted many lectures, writings, and social media posts to exactly this topic.

There are by now a great many data points bearing on various forms of the conjecture.  Unfortunately, most of them have \textbf{not} been gathered under controlled scientific conditions.  Much of the publicly available evidence is subject to severe reporting bias --- successes are announced and failures are not --- and to non-scientific incentives, with important costs and variables (the number of attempts, the amount of human scaffolding, the compute expended, the degree of contamination of the problem with prior literature) frequently left undisclosed.  Furthermore, the \emph{truth value} of a given form of the conjecture is sometimes conflated with its \emph{desirability}.  

Despite the central relevance of the AI Capability Conjecture to Question \ref{crq}, \textbf{this article is not about that conjecture}.  In this direction, I will mention only the recent results of the First Proof project \cite{firstproof-web}, an independent assessment of the capabilities of frontier AI models and harnesses on genuinely novel mathematics.  Each ``batch'' of the project consists of ten research-level problems, contributed by working mathematicians in a wide range of fields, whose solutions are known to the contributor but have never been posted anywhere online.  The second batch \cite{firstproof} was evaluated under controlled conditions against four AI systems, using models publicly accessible as of May 28, 2026, and the resulting solutions were refereed by experts for both correctness and quality of exposition.  Of the ten problems, seven received at least one passing grade --- that is, a solution judged essentially flawless or requiring only minor revisions --- from at least one system, with compute costs on the order of tens to hundreds of dollars per problem.  Further batches are planned.  

\section{The complement to the capability conjecture}

This article is instead about what one might call the ``orthogonal complement'' of the AI Capability Conjecture inside Question \ref{crq}.  To isolate that component, I will adopt the following imprecisely stated hypothesis.

\begin{hypothesis}[Working Hypothesis]\label{wh} A reasonably strong version of the AI Capability Conjecture is true: AI tools will, reasonably soon, become capable of performing a reasonable fraction of research-level mathematical tasks, with reasonable levels of success, quality, supervision, and cost.
\end{hypothesis}

The precise meaning of ``reasonable'' here is not critical for what follows.

For the remainder of this article, I ask the reader to \emph{assume} that Hypothesis \ref{wh} holds.  I am \textbf{not} asking the reader to want it to be true, to believe that it is true, or to accept it as true; what follows is a conditional analysis.  In particular, evidence for or against the Working Hypothesis is orthogonal to the discussion below.

\section{The orthogonal subquestion: our goals and values}

Once one conditions on the Working Hypothesis, a second fundamental subquestion comes into view:

\begin{question}[Goals and Values Question]\label{gvq}
What are the precise goals, objectives, and values of our mathematical community, and of the enterprise of mathematical research?  Not merely the \emph{explicit} goals that we communicate to the public, to our students, or to funding agencies, but the \emph{implicit} goals that we actually optimize for in practice?
\end{question}

In the past, we have largely delegated Question \ref{gvq} to the humanities --- to historians, philosophers, and sociologists of mathematics --- and focused our own attention on the technical content of our profession.  Assuming the Working Hypothesis, we will no longer have this luxury.  But --- as with the crisis in foundations --- I would argue that a critical examination of the question will ultimately prove highly valuable \textbf{regardless of the status of the Working Hypothesis}.  

So: what are our goals?  To my knowledge, no official list of the goals of mathematics has been systematically compiled, but here is a partial list:

\begin{itemize}
\item to solve unsolved problems, both pure and applied;
\item to develop new theories, structures, and techniques;
\item to understand the world around us;
\item to build and sustain a community of mathematicians;
\item to train the next generation of mathematicians, and to let them guide the future directions of the subject;
\item to contribute to the shared and cumulative network of mathematical knowledge;
\item to create enduring works of aesthetic value;
\item etc.
\end{itemize}

The reader is encouraged to extend this list further.

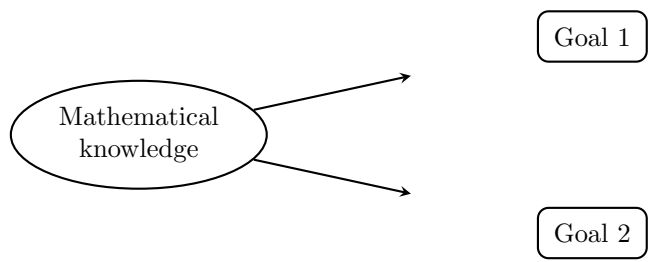
\begin{figure}[h]
\centering
\begin{tikzpicture}[>=stealth, thick]
  \node[ellipse, draw, minimum width=2.8cm, minimum height=1.4cm, align=center, font=\footnotesize] (K) at (0, 0) {Mathematical\\knowledge};
  \node[draw, rounded corners, inner sep=6pt, font=\footnotesize] (G1) at (6, 1.3) {Goal 1};
  \node[draw, rounded corners, inner sep=6pt, font=\footnotesize] (G2) at (6, -1.3) {Goal 2};
  \draw[->] (K) -- ($(K)!0.6!(G1)$);
  \draw[->] (K) -- ($(K)!0.6!(G2)$);
\end{tikzpicture}
\caption{Historically, any one goal of mathematics served as a usable proxy for the others.}\label{aligned}
\end{figure}

Historically, the above goals have been \emph{positively correlated} with one another.  Progress on any one goal has typically moved one closer to other goals as well.  For instance, in the course of solving a hard problem, a new technique may be developed, a new community of researchers forms around that technique, students are trained in it, textbooks are written, and the resulting theory turns out to be applicable elsewhere.  Because of this correlation, one could use one or two of these goals as convenient \emph{proxies} for the others, and leave the remainder implicitly stated at most.  See Figure \ref{aligned}, as well as \cite{tao-good} for a previous discussion by the author of this alignment phenomenon.

However, all metrics, when excessively optimized for, are at risk of being subjected to Goodhart's law.  In the formulation popularized by Strathern \cite{strathern}, following the original observation of Goodhart on monetary policy \cite{goodhart}:

\begin{quote}
When a measure becomes a target, it ceases to be a good measure.
\end{quote}

AI tools are particularly likely to trigger this effect, for two independent reasons.  The first is technical: generative AI is inherently \emph{ungrounded}, in the sense that it optimizes for the appearance of a satisfactory output rather than for the underlying property that the output is supposed to indicate, and so is unusually good at finding the gap between a measure and the thing it measures.  The second is economic: the financial incentives of the AI industry reward demonstrable, quotable, benchmarkable achievement on precisely the sort of metrics we have historically used as proxies.

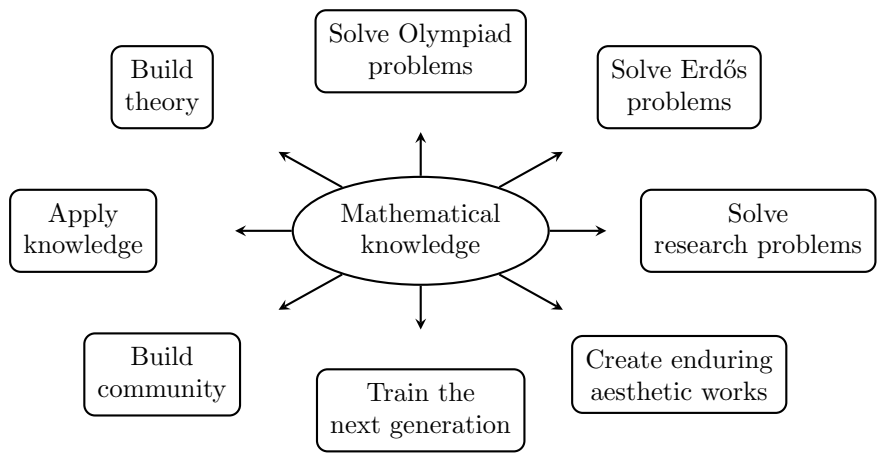
\begin{figure}[h]
\centering
\begin{tikzpicture}[>=stealth, thick, scale=0.95]
  \node[ellipse, draw, minimum width=2.6cm, minimum height=1.3cm, align=center, font=\footnotesize] (K) at (0, 0) {Mathematical\\knowledge};
  \node[draw, rounded corners, inner sep=5pt, align=center, font=\footnotesize] (G1) at (4.7, 0)     {Solve\\research problems};
  \node[draw, rounded corners, inner sep=5pt, align=center, font=\footnotesize] (G2) at (3.6, 2.0)   {Solve Erd\H{o}s\\problems};
  \node[draw, rounded corners, inner sep=5pt, align=center, font=\footnotesize] (G3) at (0, 2.5)     {Solve Olympiad\\problems};
  \node[draw, rounded corners, inner sep=5pt, align=center, font=\footnotesize] (G4) at (-3.6, 2.0)  {Build\\theory};
  \node[draw, rounded corners, inner sep=5pt, align=center, font=\footnotesize] (G5) at (-4.7, 0)    {Apply\\knowledge};
  \node[draw, rounded corners, inner sep=5pt, align=center, font=\footnotesize] (G6) at (-3.6, -2.0) {Build\\community};
  \node[draw, rounded corners, inner sep=5pt, align=center, font=\footnotesize] (G7) at (0, -2.5)    {Train the\\next generation};
  \node[draw, rounded corners, inner sep=5pt, align=center, font=\footnotesize] (G8) at (3.6, -2.0)  {Create enduring\\aesthetic works};
  \foreach \g in {G1,G2,G3,G4,G5,G6,G7,G8}{\draw[->] (K) -- ($(K)!0.55!(\g)$);}
\end{tikzpicture}
\caption{Under excessive optimization, the goals of mathematics diverge from each other.  The diagram is of course extremely oversimplified; in particular, it should be very much higher dimensional.}\label{diverge}
\end{figure}

Consequently, excessive optimization for one or two goals may cause the many previously aligned goals of mathematics to \emph{diverge} from one another; see Figure \ref{diverge}.  

\section{A case study: problem solving}

To make the preceding discussion concrete, I will focus on a single component of mathematical research: \emph{problem solving}.  It should be stressed that this is not at all the only aspect of our profession.  \emph{Theory building}, for instance, is a complementary activity of at least equal significance, and one that requires its own separate analysis; so do teaching, mentoring, and the many forms of service by which a research community sustains itself.  But problem solving is a natural first case study, both because it is the aspect most susceptible to being impacted under the Working Hypothesis, and because it is the aspect for which our implicit goals are the furthest from our explicit ones.

Suppose then that we try to write down what we want from problem solving.  A first attempt might be the following.

\begin{goal}[first attempt]\label{goal1}
Solve as many unsolved problems as possible.
\end{goal}

Under Goal \ref{goal1}, we are trying to optimize the flow in a very simple network:
\[
\begin{tikzcd}[column sep=huge]
  \boxed{\text{Open problems}} \arrow[r, "\text{proof generation}"] & \boxed{\text{Solutions}}
\end{tikzcd}
\]
Even before the advent of AI, we knew that this metric was inadequate, and we knew it for an entirely mundane reason: optimizing it produces a large number of \emph{incorrect} solutions to major open problems.  Every working mathematician with a public email address is familiar with the steady stream of purported proofs of the Riemann hypothesis.  Hence, we may update our goal:

\begin{goal}[second attempt]\label{goal2}
Solve as many unsolved problems as possible, and \emph{verify them to be correct}.
\end{goal}

The corresponding network acquires a second step:
\[
\begin{tikzcd}[column sep=huge, row sep=large]
  \boxed{\text{Open problems}} \arrow[r, "\text{proof generation}"] & \boxed{\text{Unverified solutions}} \arrow[d, "\text{proof verification}"] \\
  & \boxed{\text{Verified solutions}}
\end{tikzcd}
\]
Advances in AI, and in autoformalization into proof assistant languages\footnote{For further discussion of recent developments in formalization, see \cite{avigad}.} such as Rocq, HOL, or Lean \cite{lean, mathlib}, have significantly accelerated both proof generation and proof verification in many cases, and under the Working Hypothesis this acceleration will continue.  A formally verified proof is, after all, precisely a proof whose correctness no longer depends on the reputation or the diligence of its author.

But now a new failure mode appears.  What if an AI tool generates a lengthy proof that is verified to be correct, but which nobody --- not even the humans who prompted the tool --- understands?  This is no longer hypothetical.  Sites devoted to collecting mathematical problems, such as the Erd\H{o}s problems database \cite{erdosproblems},\footnote{For a study of the emergence of flourishing online mathematical communities, see \cite{pease}.} already contain dozens of AI-generated proof submissions.  Many of these are likely to be correct; but in a substantial number of cases no human expert has yet volunteered to verify and vouch for them, and in several cases the human submitters have themselves declared that they are not qualified to do so.  We may soon be faced with the very real possibility of a verified proof of a major result that no human understands well enough to explain.

Thus, we may update our goal again:

\begin{goal}[third attempt]\label{goal3}
Solve as many unsolved problems as possible, verify them to be correct, and \emph{ensure that the results can be clearly communicated to and understood by the mathematical community}.
\end{goal}

\[
\begin{tikzcd}[column sep=huge, row sep=large]
  \boxed{\text{Open problems}} \arrow[r, "\text{proof generation}"] & \boxed{\text{Unverified solutions}} \arrow[d, "\text{proof verification}"] \\
  \boxed{\text{Well-written solutions}} & \boxed{\text{Verified solutions}} \arrow[l, "\text{proof exposition}"']
\end{tikzcd}
\]

Current AI tools have a decidedly mixed record with proof exposition.  On the one hand, the spelling, the grammar, and the formatting are close to flawless\footnote{One can argue that they are \emph{too} flawless.}.  On the other hand, the writing very often dwells at length on trivialities while passing briefly through --- or even actively obscuring --- the most interesting and novel portions of the argument.  AI-generated mathematical texts also frequently fail to situate the result in the prior literature, or to offer the high-level overview that lets a reader decide whether the argument is worth their time.

Proof exposition is admittedly a much ``fuzzier'' optimization target than proof verification, and the Working Hypothesis predicts that AI tools will improve at it considerably from current levels.  But here I want to make a point that I think is under-appreciated: exposition, too, can be over-optimized.  A proof can be \textbf{too slickly written}, with the routine steps and the genuinely difficult steps presented as being equally easy to digest.

\begin{figure}[h]
\centering
\includegraphics[trim=0 200 0 0, clip, width=0.75\textwidth]{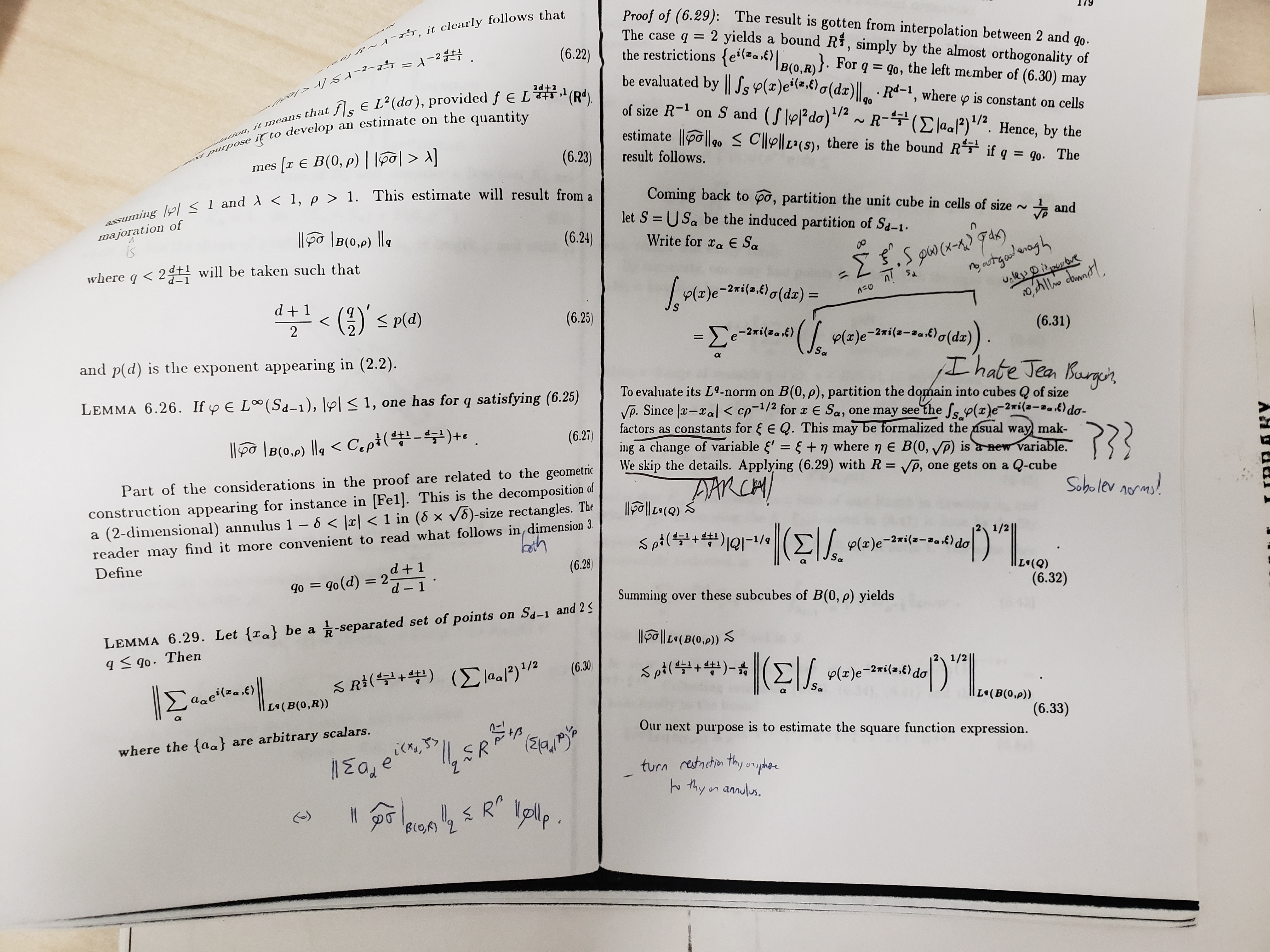}
\caption{A page from a 1991 paper of Bourgain \cite{bourgain}, annotated by my much younger (and very frustrated) self.  But by fighting my way through these texts, I came to understand Bourgain's way of thinking, and in time I actively sought out his papers to read.  See also \cite{bourgain-tricks}, \cite{bourgain-obit}.}\label{annotated}
\end{figure}

In a human-written proof, the parts of the argument that the author found difficult typically retain some \emph{natural friction}: an apologetic remark, an unusually careful lemma, a change of notation, a paragraph that has clearly been rewritten several times.  This friction is informative.  It signals to the reader where to slow down and pay attention, and it is one of the main channels by which the tacit knowledge of a field is transmitted.  An excessively AI-polished proof may sand away both the ``artificial'' friction (typos, awkward phrasing, disorganization) and the ``natural'' friction, leaving a text that is easy to read and hard to learn from.  Paradoxically, the ``mistakes'' in human exposition can be genuinely helpful to the reader; see Figure \ref{annotated}.

It is worth recalling Thurston's formulation of the point, from his classic essay \cite{thurston}, which remains as relevant in the age of AI as it did in 1994:

\begin{quote}
``We are not trying to meet some abstract production quota of definitions, theorems and proofs.  The measure of our success is whether what we do enables people to understand and think more clearly and effectively about math.''
\end{quote}

For a proof to actually contribute to its field, then, it is not enough for it to be correct, and not enough for it to be readable.  It also needs to be \emph{accepted} and \emph{valued} by the community: other mathematicians need to digest the result and incorporate it into their own work.  Authors can materially assist in this digestion process, by describing the insights, the false starts, and the stories from the period when they were working on the problem.  In contrast, current AI tools are quite opaque about their own problem-solving process, and this is particularly true of proprietary models whose inner workings are a corporate secret.  

Thus, we may update our goal yet again:

\begin{goal}[fourth attempt]\label{goal4}
Solve unsolved problems, verify them to be correct, ensure they are clearly communicated, and have them \emph{digested and accepted by the mathematical community}.
\end{goal}

\[
\begin{tikzcd}[column sep=huge, row sep=large]
  \boxed{\text{Open problems}} \arrow[r, "\text{proof generation}"] & \boxed{\text{Unverified solutions}} \arrow[d, "\text{proof verification}"] \\
  \boxed{\text{Well-written solutions}} \arrow[d, "\text{proof publication}"'] & \boxed{\text{Verified solutions}} \arrow[l, "\text{proof exposition}"'] \\
  \boxed{\text{Accepted solutions}} &
\end{tikzcd}
\]

Community acceptance of a result is, by its nature, slow and human.  It can be \emph{encouraged} by good exposition and careful writing, but it is ultimately an external process that \textbf{cannot be optimized purely by the authors and their tools}.  Our current publication infrastructure relies on human editors and referees to provide this acceptance, voluntarily and largely without credit.  This work is routinely regarded as less prestigious than the work of generating proofs in the first place; but it is an essential component of the profession, and it is precisely the mechanism by which the individual achievements of mathematicians are converted into collective progress and understanding.

AI evaluation tools may well serve as useful \emph{filters} in this process --- one can imagine journals automatically triaging submissions that are flagged for inadequate verification, missing attribution, or incoherent exposition, in the same way that plagiarism detection is used today.  Such filters, while controversial to implement, would conserve the scarce resource of expert human attention.  But passing an automatic filter is not a substitute for community acceptance; I do not believe that human referees can be removed from the publication process.

Finally, even publication is not the last stage.  Key results should ultimately become part of the definitive textbooks and reference material of their subject, in the form in which they are taught to the next generation of students.  This process of \emph{canonicalization} --- in which a result is restated in its natural generality, given its right proof rather than its first proof, connected to its neighbors, and absorbed into the standard toolkit --- is the slowest stage of all.  It requires broad, deliberative consensus, and it is the stage \emph{least} amenable to optimization by AI tools.  It is also, in my view, the \emph{most valuable} part of the entire process.  Many applications of mathematics only become feasible once the underlying theory has been fully digested in this way.  Indeed, the very success of AI tools in mathematics depends crucially on the canonical theories that human mathematicians have painstakingly built and rebuilt over the centuries: the training data for these tools is, quite literally, the output of the canonicalization process.

This gives us a (potentially) final version of the problem solving goal:

\begin{goal}[final attempt?]\label{goal5}
Solve unsolved problems, verify them to be correct, ensure they are clearly communicated, and have them digested, accepted, and \emph{incorporated into the definitive theory of the field}.
\end{goal}

\begin{figure}[h]
\centering
\begin{tikzcd}[column sep=huge, row sep=large]
  \boxed{\text{Open problems}} \arrow[r, "\text{proof generation}"] & \boxed{\text{Unverified solutions}} \arrow[d, "\text{proof verification}"] \\
  \boxed{\text{Well-written solutions}} \arrow[d, "\text{proof publication}"'] & \boxed{\text{Verified solutions}} \arrow[l, "\text{proof exposition}"'] \arrow[d, bend right=85, looseness=6, dotted, "\text{proof digestion}"] \\
  \boxed{\text{Accepted solutions}} \arrow[r, "\substack{\text{proof} \\ \text{canonicalization}}"'] & \boxed{\text{Definitive solutions}}
\end{tikzcd}
\caption{The problem-solving pipeline, in the form arrived at by iterating Goals \ref{goal1}--\ref{goal5}.  What began as a single arrow has become a chain of five stages, of which only the first was ever an explicit goal of the community.}\label{pipeline}
\end{figure}
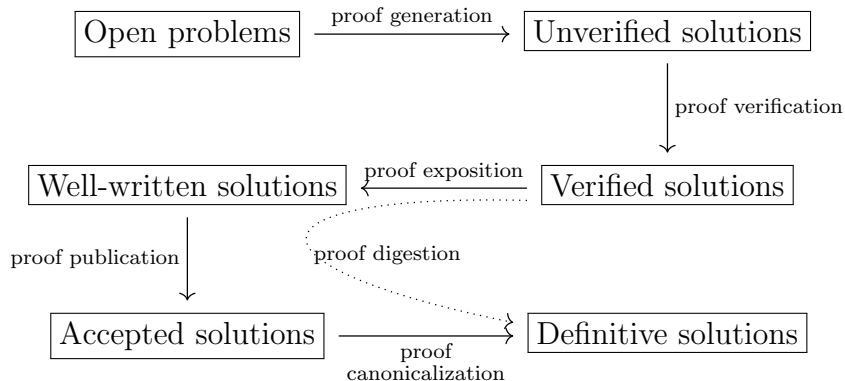

The specific pipeline in Figure \ref{pipeline} may still be oversimplified and subject to further analysis; however it illustrates the nuances one uncovers when one deconstructs a goal, such as problem solving, which seems simple on the surface, but in fact carries many implicit subgoals that are worth making explicit.

\section{Proof scarcity and proof abundance}

If the Working Hypothesis holds, then in the absence of suitable policy and cultural changes, significant ``impedance mismatches'' --- or, to use a less flattering metaphor, \emph{proof indigestion} --- will emerge all along the pipeline of Figure \ref{pipeline}:

\begin{itemize}
\item AI-generated proofs will accumulate faster than they can be verified;
\item verified AI-generated proofs will accumulate faster than they can be given a readable write-up;
\item AI-generated proofs, even those required to be both correct and well written, will overwhelm a traditional peer review system that depends on volunteer expert labor;
\item and even the published proofs will be too numerous for the community to work into definitive form.
\end{itemize}

In short, we will transition from an era of \emph{proof scarcity} to an era of \emph{proof abundance}.  Most of our institutions --- journals, priority conventions, hiring and promotion criteria, prizes, the very notion of a research program --- were designed under the assumption of scarcity, and it should not surprise us if they behave poorly under abundance.  Some signs of this indigestion were already appearing before the advent of modern AI: the growth in the volume of the literature, the increasing length and specialization of major proofs, and the well-documented strain on the refereeing system all predate the present moment.  But the advent of AI will exacerbate these existing stresses markedly.

\section{From goals to recommendations}

Identifying the goals of problem solving in the way we have done above makes it considerably easier to see how to respond to these emerging impedance mismatches, because each mismatch is now attached to an identified stage and an identified value.  I do not intend to propose a full program here.  Instead, I will point to the Leiden Declaration on Artificial Intelligence and Mathematics \cite{leiden}, published in June 2026 and endorsed by the International Mathematical Union, which I regard as an excellent starting point.  The declaration arose from a 2025 workshop at the Lorentz Center in Leiden, and consists of twenty-three recommendations addressed to individual mathematicians, to mathematical organizations and not-for-profit funders, and to policymakers.  

Rather than reproduce the declaration, let me quote four of its recommendations to individual mathematicians, and offer a commentary on each from the perspective developed above.

\medskip

\begin{quote}
\textbf{Disclose tool use.} Transparently disclose the use of automated tools, including large language models, machine learning systems, proof assistants, and other mathematical software. Include a ``Tool and computational resource disclosure'' section in your papers; many journals, publishers, and professional organizations have already developed guidelines for this, and though the precise form of such a section will necessarily evolve, we encourage authors to live up to the spirit reflected in the UNESCO Recommendation on Open Science \cite{UNESCOOpenScience} and the FAIR principles \cite{WilkinsonEtAlFAIR}. When acting as a reviewer, abide by publisher guidelines. If the use of artificial intelligence is allowed, be transparent about how you used it, and take responsibility for any significant recommendations you make.
\end{quote}

\medskip

The scenario to be avoided at all costs is one in which authors use AI tools covertly to aid their work, but conceal that usage in order to avoid criticism from their peers.  (For my own disclosure of AI tools in preparing this paper, see Section \ref{ack}.)

\medskip

\begin{quote}
  \textbf{Support the needs of reviewing.} The use of artificial intelligence in preparing papers can introduce material that makes reviewing more demanding. Make it easier for your peers to review your work by disclosing tool use, giving precise and complete references to previous results, and providing formal proofs where feasible and appropriate.
\end{quote}

\medskip

More broadly, I argue that we need to \emph{decrease} the emphasis that our culture places on proof generation, and in particular on being the ``first'' to solve a problem, and correspondingly \emph{increase} the emphasis we place on proof digestion: exposition, refereeing, publication, and canonicalization.  See also the recent essay of Bessis \cite{bessis} on the need to move away from the ``theorem economy'' based primarily on proof generation. 

\medskip

\begin{quote}
  \textbf{Affirm the humanity of authorship.} Credit and responsibility continue to belong to humans within the mathematical community and should not be given to automated systems. Artificial intelligence may obscure, but does not replace, the collective human labor behind a result.
\end{quote}

\medskip

\begin{quote}
\textbf{Put effort into proper attribution.} The known limitations of automated tools in properly attributing ideas create a corresponding obligation for proactive effort to find and credit the sources that made a new result possible. Where a satisfactory attribution is not possible, state this explicitly in the publication.
\end{quote}

\medskip

My own suggested rule of thumb: \emph{if the authors cannot convincingly demonstrate that they are able to give a clear, expert-level talk on their results, one that is correct and properly attributed, then the result should not be published}.  A proof that no human can properly explain should be viewed as incomplete, even if it has been formally verified.

\section{Closing thoughts}

I have presented problem solving as \emph{one} aspect of mathematics in which the Working Hypothesis forces us to inspect goals and values that we have long been able to leave implicit.  But the Working Hypothesis potentially impacts many other aspects of our work --- teaching, mentoring, hiring, grant applications, refereeing, public outreach --- and a similar analysis should be performed for each of them.  

The conclusions of such analyses will not be uniform.  In some areas, particularly in education and in the training of young mathematicians, it will be crucial to emphasize the irreducibly human aspect of our work, and to restrict the use of AI tools quite tightly; the goal of training a mathematician is not achieved by producing correct homework.  In other areas, we will need to take the initiative on AI usage, and define best practices for incorporating these tools into our workflows on our own terms rather than on terms set for us by vendors.  We will also need new workflows and new infrastructures to complement our traditional ones --- collaborative formalization projects, structured problem databases, new venues for exposition and for the publication of negative or partial results; see Appendix \ref{app} for a partial list.

Above all, our community needs to come together to have open and honest discussions about \emph{both} of the subquestions identified here: about AI capability, and about our own goals and values.  This is again a recommendation of the Leiden declaration:

\medskip

\begin{quote}
\textbf{Participate in public discourse.} Mathematicians have a responsibility to support serious science journalism and to engage in public discourse to explain and contextualize artificial intelligence-assisted methods and results. This is particularly important for work within our own subfields, where specialized knowledge is required to assess claims about the depth, difficulty, and significance of results. Moreover, we encourage mathematicians to seek opportunities to cooperate with and support other researchers and creative professionals facing similar challenges.
\end{quote}

\section{Acknowledgments}\label{ack}

This article is based on a public lecture delivered at the International Congress of Mathematicians in July 2026.  I thank Bryna Kra, Jeremy Avigad, Martin Hairer, Akshay Venkatesh, and Emily Riehl for their feedback on early versions of that lecture.

AI assistance was used to perform literature search, to generate diagrams, to autocomplete text, and to convert the slides into a paper format.

\appendix

\section{Some new workflows and infrastructures}\label{app}

For the interested reader, I list a few existing projects that illustrate the kinds of new infrastructure discussed above:

\begin{itemize}
\item Mathlib, a unified formalized library of mathematics: \url{https://mathlib.org/} \cite{mathlib};
\item Mathematical Discourse, a peer-reviewed video journal for research talks: \url{https://www.mathematicaldiscourse.org/} \cite{discourse};
\item the Erd\H{o}s problems database: \url{https://www.erdosproblems.com} \cite{erdosproblems};
\item the optimization constants database:\\ \url{https://github.com/teorth/optimizationproblems} \cite{optimization};
\item the SAIR Foundation mathematics competitions:\\ \url{https://competition.sair.foundation/competitions} \cite{sair};
\item the First Proof project: \url{https://1stproof.org/} \cite{firstproof-web}.
\item the Palomar formalized proof registry: \url{https://palomar-registry.org/}
\end{itemize}

\bibliographystyle{amsplain}

\end{document}